\documentclass[12pt]{paper}
 \usepackage{epsfig}
\usepackage{graphics}
 \usepackage{color}

\newcommand{\Mo}{\mathbb{M}}

\usepackage{color}
\usepackage{epsfig,cmmib57}
\usepackage{amssymb,amsmath,exscale}
 \numberwithin{equation}{section}
\newcommand{\ZLL}{\langle\!\langle}
\newcommand{\ZRR}{\rangle\!\rangle}

\newcommand{\ZEP}{\epsilon}

\newcommand{\ZSUno}{\sum _{n=1}^{+\ZIN}}

\newcommand{\ZOMq}{\Omega}

\newcommand{\ZT}{{\mathcal T}}
\newcommand{\ZA}{{\mathcal A}}

\newcommand{\zg}{\gamma}

\newcommand{\intT}{\int_0^T}
\newcommand{\intt}{\int_0^t}
\newcommand{\ints}{\int_0^s}

\newtheorem{Theorem}{Theorem}

\newtheorem{Lemma}[Theorem]{Lemma}

\newtheorem{Remark}[Theorem]{Remark}

\newcommand{\zdiaform}{\mbox{~~\zdia}}

\newcommand{\zaa}{\alpha}

\newcommand{\zt}{\tau}
\newcommand{\zdia}{~~\rule{1mm}{2mm}\par\medskip}

\newcommand{\ZLA}{\label}
\newcommand{\ZIN}{\infty}
\newcommand{\zProof}{{\noindent\bf\underbar{Proof}.}\ }

\newcommand{\zzr}{{\rm I\hskip-2.1pt R}}
\newcommand{\ZBI}{\bibitem}
\newcommand{\ZD}{\;\mbox{\rm d}}
\newcommand{\zl}{\lambda}

\author{
L. Pandolfi\thanks{Dipartimento di Scienze Matematiche ``Giuseppe Luigi Lagrange'', Politecnico di Torino, Corso Duca degli Abruzzi 24, 10129 Torino, Italy (luciano.pandolfi@polito.it)}
}
\title{Controllability of a viscoelastic plate using one boundary control in displacement or bending\thanks{
This papers fits into the research program of the GNAMPA-INDAM and has been written in the framework of the   ``Groupement de Recherche en Contr\^ole des EDP entre la France et l'Italie (CONEDP-CNRS)''.}}
\begin{document}
 
 \maketitle 
 
\noindent {\bf\underline{Abstract}:}  In this paper we consider a viscoelastic plate (linear viscoelasticity of the Maxwell-Boltzmann type) and we compare its controllability properties with the (known) controllability of a purely elastic plate (the control acts on the boundary displacement or bending). By combining operator and moment methods, we prove that the viscoelastic plate inherits the   controllability properties of the purely elastic plate.
 \medskip
 
\section{Introduction} 

Controllability of elastic and viscoelastic bodies is a standing subject of investigation in systems theory and in particular controllability of an elastic plate whose (vertical) displacement is described by
\begin{equation}
\ZLA{eq:sistemaelastico}
u''+\Delta^2 u=F
\end{equation}
(and   controls acting on different boundary conditions) has been studied in many papers  after the first results in~\cite{lionsLIBRO,LionsNeuman}.

In Eq.~(\ref{eq:sistemaelastico}), the apex denotes time derivative, $\Delta$ is the laplacian,   $u=u(x,t)$ denotes the vertical displacement at time $t$ and position $x$ of the plate. So $x\in \ZOMq$, and ${\rm dim}\,\ZOMq=2$ in the physically significant cases. 

We assume that $\ZOMq$ has $C^3$ boundary and   we associate the following initial and boundary conditions to Eq.~(\ref{eq:sistemaelastico}):
\begin{equation}
\ZLA{eq:elasticoBoundCondi}
\begin{array}{l}
u(0)=u_0\,,\ u'(0)=u_1\, \quad {\rm and}\quad\\[2mm] 
 \left\{\begin{array}{ll}
  \mbox{either {\bf case (A)}:}&   \zg_0u=g\,,\ \zg_1u=0  \\
\mbox{or {\bf case (B)}:}&  \zg_0u=0\,,\ \zg_1u=g\,.  
\end{array}\right.
\end{array}
\end{equation}
The symbols $\zg_0$ and $\zg_1$ denote the traces on $\partial\ZOMq$ of $u$ and of its normal derivative.

 It is well known that system~(\ref{eq:sistemaelastico}) is controllable (in a suitable space $X$ described below) using  square integrable controls $g $, i.e., for every $(u_0,u_1)\in X$ and every target $(\hat u_0,\hat u_1)\in X$ there exists a control $g  $ such that $(u(T),u'(T))=(\hat u_0,\hat u_1)$ (note that in the study of controllability we can use $F=0$, $u_0=0$, $u_1=0$). 

The noticeable fact is that $T>0$ is arbitrary, as first proved in~\cite{ZuazuaCR}.

Important references for the previous result are~\cite{Komornik,LasTrigg}. See also~\cite{AASSILA,Arena,Bayili,Guo,KUK,Niane,Shubov} for extensions and the case that the control acts in different boundary conditions. We note that  controllability in {\bf case (A)} is studied in~\cite{LasTrigg} when $\Gamma=\partial\ZOMq$ (see also~\cite{Niane})  while {\bf case (B)} is studied in~\cite{Komornik} when $\Gamma$ is a suitable subset of $\partial\ZOMq$. So, our standing assumption is that \emph{the controls $g$ are supported in  $\Gamma\subseteq \partial\ZOMq$ and that the system~(\ref{eq:sistemaelastico})-(\ref{eq:elasticoBoundCondi}) is controllable.}

Our goal   is the proof that controllability of the purely elastic plate is inherited by a viscoelastic plate whose dynamic is described by the equation
\begin{equation}\ZLA{eq:sistemsvisco}
\left.\begin{array}{l}w''+\Delta^2w+\intt M(t-s)\Delta^2 w(s)\ZD s=F\,,\\
w(0)=w_0\,, \ w'(0)=w_1\ \mbox{and} \\
 \left\{\begin{array}{ll}
  \mbox{either {\bf case (A)}:}&   \zg_0w=g\,,\ \zg_1w=0  \\
\mbox{or {\bf case (B)}:}&  \zg_0w=0\,,\ \zg_1w=g\,.  
\end{array}\right.
\end{array}
\right.
\end{equation}
We are going to prove that controllability holds in the same space $X$   and at every  time $T>0$, as for the elastic plate.

The plan of the paper is as follows: 
previous references and 
needed results on the Eq.~(\ref{eq:sistemaelastico}) are commented both in this introduction and in the next section, where we introduce preliminary notations,  assumptions    and   the spaces $X$ in which controllability is studied.

Section~\ref{sect:TheSOlutions} studies the solution of system~(\ref{eq:sistemsvisco}) while controllability is studied in Sect.~\ref{sect:ControllabilityVISCO}.

The control problems {\bf (A)} and {\bf (B)} are studied in parallel, with similar methods and results (but the control spaces are different).  

 \begin{Remark}
 {\rm
 We note:
 \begin{itemize}
 \item we shall use condition~(\ref{eq:stimaautov}) below, which holds in the physical case   ${\rm dim}\,\ZOMq =2$ but arguments similar to those in~\cite{PandLIBRO,PandLame} can be used to remove this condition.
 \item we expect that similar arguments as in this paper can be used to extend to the viscoelastic case the existing controllability results for an elastic plate, under  controls acting on different boundary conditions, for example in the boundary moment. This will be the subject of future investigations. 
 \item an active field of research now is controllability of the connection of distributed   systems, see~\cite{AASSILA} for the study of controllability of two connected plates. It seems that this kind of problems has not been considered in the case of viscoelasticity.
 \end{itemize}
  
 }
 \end{Remark}
 
We finish this section by mentioning   few of the existing results on controllability of viscoelastic plates. The first result seems to be in~\cite{LeugBEAM} (a paper which studies controllability of a beam). In this paper the control acts on the moment (i.e. $w_{xx}(0,t)$ is controlled, while it is assumed $w(0,t)=w(L,t)=w_{xx}(L,t)=0$) and controllability is proved in the space $\left (H^2(0,L)\cap H^1_0(0,L)\right )\times L^2(0,L)$. This same problem has been studied in~\cite{AvdoninBELI} in the space $H^1_0(0,L)\times H^{-1}_0(0,L)$.
The result in~\cite{LeugBEAM} has then been extended to a rectangle in~\cite{LeugRectangle}.
As we shall see, our contollability result for the viscoelastic plate is obtained from the corresponding result of the elastic plate via perturbation. This idea was first used in~\cite{LagLIONs,LasVisco} by assuming  that the memory kernel is (smooth and) sufficiently small (while we shall use a compactness argument).   The ``smallness'' assumption  was removed in~\cite{KimPlate} (in {\bf case (B)}) but the methods in this paper require the restrictive assumption that the control acts on the entire $\partial\ZOMq$ (a condition not required to control the elastic plate, see~\cite{Komornik}, and not required by   the result we prove in this paper). The ideas that we use here are very different from those in~\cite{KimPlate} and relay on a combination of operator and moment methods.

\section{Notations, assumptions and preliminaries}

As already stated, we assume that $\ZOMq\subseteq\zzr^2$ is a bounded region with $C^3$ boundary and we assume that the memory kernel $M(t)$, defined on $[0,T]$, is of class $C^2$.
The regularity of $\partial\ZOMq$ is used to ensure that every solution of~(\ref{eq:sistemaelastico}) or~(\ref{eq:sistemsvisco}) is the limit of smooth solutions, so that it is possible to compute with smooth fuctions and pass to the limit, for example in the proofs of the direct and inverse inequalities introduced below.

We shall use consistently the following notations: 

\begin{itemize}
\item $u$ or $\phi$   denote  the solution of the purely elastic system~(\ref{eq:sistemaelastico}). The greek letter $\phi$ is used when we want to stress that the control is put equal to zero, $g=0$. If we need to stress the dependence of $u$ on the control $g$ then we write $u_g$.
\item the solution of the viscoelastic system~(\ref{eq:sistemsvisco}) is consistently denoted $w$ or $\psi$ (we use $\psi$ when $g=0$) and we use $w_g$ to denote the dependence of $w$ on the boundary control $g$.
\item We use the operator $A=\Delta^2$ in $L^2(\ZOMq)$ with
\[
{\rm dom}\, A=\left \{ \phi\in H^4(\ZOMq)\,,\quad \zg_0\phi=0\,,\ \zg_1\phi=0\right \}\,.
\] 
The operator $A$ is selfadjoint positive with compact resolvent so that there exists an orthonormal basis $\{\phi_n\}$ of eigenvectors of $A$
\[
A\phi_n=-\zl_n^2\phi_n\,.
\]
When the eigenvectors are ordered in a nondecreasing sequence then there exists $\zaa>0$ such that  $\zl_n^2\geq \zaa n^{4/d}$, $d={\rm dim}\,\ZOMq$  (see~\cite[Sect.s~13-14]{Agmon}) so that in our case $d=2$ we have
 \begin{equation}
 \ZLA{eq:stimaautov}
 \zl_n^2\geq \zaa n^2\,,\qquad \ZSUno \frac{1}{\zl_n^2}<+\ZIN\,.
 \end{equation}
Finally, we introduce
\[
\ZA=A^{1/2}\,,\qquad R_+(t)=\frac{1}{2}\left [ e^{\ZA t}+e^{-\ZA t}\right ]\,,
\quad 
R_-(t)=\frac{1}{2}\left [ e^{\ZA t}-e^{-\ZA t}\right ]\,.
\]
\end{itemize}

We use $\Gamma$ to denote $\partial\ZOMq$   or the relatively open subset of $\partial\ZOMq$ on which the control acts (in its relative complement boundary conditions will be homogeneous) and for every fixed $T>0$ we put $\Sigma=\Gamma\times (0,T)$. We denote $\ZD\Gamma$ and $\ZD\Sigma$ the surface measures.
 
 \subsection{The space of controllability}
 Controllability in the cases {\bf (A)} and {\bf(B)} can be studied with parallel arguments, but the control space $X$ is not the same. We introduce the following spaces (note that ${\rm dom}\, A^{1/2}=H^2_0(\ZOMq)$,  ${\rm dom}\, A^{1/4}=H^1_0(\ZOMq)$):
\[
\begin{array}
{ll}
\mbox{\bf case   (A)}& 
 Y=     \left ({\rm dom}\, A^{3/4}\right )\times \left ({\rm dom}\,A^{1/4}\right )\,,\quad 
 \tilde Y=\left ({\rm dom}\,A^{1/4}\right ) \times \left ({\rm dom}\, A^{3/4}\right )\,, 
 \\
&X=\tilde Y' =\left ({\rm dom}\,A^{1/4}\right )'\times \left ({\rm dom}\, A^{3/4}\right )'\,,\\[3mm]
\mbox{\bf  case  (B)}&
Y
= {\rm dom }\,A^{1/2}\times L^2(\ZOMq)\,,\quad  
 \tilde Y=
L^2(\ZOMq)\times {\rm dom}\,A^{1/2}\,,
 \\ & X 
=\tilde Y'
 =L^2(\ZOMq)\times \left ({\rm dom}\,A^{1/2}\right )' \,.
\end{array}
\]

 In both the cases, the following result holds (see~\cite{LasTrigg}  for the case {\bf (A)} and~\cite{Komornik} for the case~{\bf(B)}):
 \begin{Theorem}\ZLA{teo:proprieREGOela}
We have:
\begin{itemize}
\item if $g\in\mathcal{D}(\Sigma)$, $F\in \mathcal{D}(\ZOMq\times(0,T))$,  $(u_0, u_1)\in Y$ then $(u,u')\in C ([0,T], Y)$. 
\item if $(u_0,u_1,F,g)\in X\times L^2(\ZOMq\times(0,T))\times L^2(\Sigma)$ then $(u,u')\in C\left ([0,T],X \right )$ and $(u_0,u_1,F,g)\mapsto (u,u')$  is continuous  in the indicated spaces.
 \item Let
 \[
\begin{array}
{ll}
\mbox{\bf case~(A)} &  \ZT\phi=-\zg_1\Delta\phi\\
\mbox{\bf case~(B)} &\ZT\phi=\zg_0\Delta \phi\,.
\end{array} 
 \]
 Then, for every $T>0$ there exists $M=M_T$ such that the  following \emph{direct inequality} holds:
\begin{equation}\ZLA{eq:direct}
\int_\Sigma \left |\ZT\phi\right |^2\ZD\Sigma\leq M\left (
\|(\phi_0,\phi_1)\|^2_Y+\|F\|^2_{L^1(0,T;L^2(\ZOMq))}
\right )   \,.
\end{equation}
\end{itemize}
\end{Theorem}

We refer to~\cite{Komornik,LasTrigg,LasTriggBUMI} for these results and for the results below.

In the study of controllability we can assume $u_0=0$, $u_1=0$, $F=0$.
Let us assume that the support of the control $g$ be contained in a (relatively open) subset $\Gamma\subseteq\partial\ZOMq$. Controllability is the property that the map $\Lambda_T $ from $L^2(0,T;L^2(\Gamma))$ to $X$:
\[
g\mapsto\Lambda_T g= (u_g(T),u'_g(T))
\]
 is surjective and this is equivalent to the property that $\Lambda_T^*$ is coercive, i.e. that the following      \emph{inverse inequality} holds:
\begin{equation}\ZLA{eq:inverse}
m\|(\phi_0,\phi_1)\|^2_Y 
 \leq \int_\Sigma \left |\ZT\phi\right |^2\ZD\Sigma\,,\qquad m>0\,.
\end{equation} 
This inequality  in the case~{\bf(A)} is proved in~\cite{LasTrigg} when $\Gamma=\partial\ZOMq$
 (see also~\cite{Bayili})  and in   case~{\bf(B)} it is proved in~\cite{Komornik} ($\Gamma$ is a suitable part of $\partial\ZOMq$).

\emph{In essence, our goal  in this paper is the proof that both the direct and inverse inequalities extend to the system~(\ref{eq:sistemsvisco}).}
 
\subsection{The solutions of the elastic system}

 The solutions of Eq.~(\ref{eq:sistemaelastico}) with conditions~(\ref{eq:elasticoBoundCondi}) are given by
\begin{align}
\nonumber
u(t)&=R_+(t) u_0+\ZA^{-1}R_-(t) u_1+\ZA^{-1}\intt R_-(t-s) F(s)\ZD s-\\
\ZLA{eq:soluELASTCO}&-\ZA\intt R_-(t-s) Dg(s)\ZD s
\end{align}
where $D$: $L^2(\partial\ZOMq)\mapsto L^2(\ZOMq)$   is defined by
\[
u=Dg\ \iff \Delta^2 u=0 \ {\rm and}\ \left\{\begin{array}{ll}
\mbox{\bf case (A)} & \zg_0u=g\,, \ \zg_1 u=0\\
\mbox{\bf case (B)} &\zg_0 u=0\,, \zg_1 u=g
\,.
\end{array}\right.
\]
It is known that $D$ is continuous and in fact even compact: it takes values in   ${\rm dom }\,A^{1/8}$ in the case {\bf (A)} and in ${\rm dom}\,A^{3/8+\ZEP}$ (any $\ZEP>0$) in the case {\bf (B)}.

%

%

The following equality is known both in case {\bf(A)} and in case~{\bf (B)} (see~\cite{LasTrigg}):
\begin{equation}\ZLA{eq:operatorefrontiera}
\mbox{if $\phi\in {\rm dom}\, A$ then $-D^*A\phi=\ZT\phi$}\,.
\end{equation}
Observe that $(\phi_0,\phi_1)\mapsto(\phi(t),\phi'(t))$ (solution of~(\ref{eq:sistemaelastico}) with zero boundary condition and $F=0$) is a $C_0$-semigroup on ${\rm dom }\,A^{1/2}\times L^2(\ZOMq)$ and so if $(\phi_0,\phi_1)\in \mathcal{D}\times\mathcal{D} $ then $\phi(t)\in {\rm dom}\,A$ for every $t$ and we have $-D^*A\phi(t)=\ZT\phi(t)$.  So, the direct inequality in particular implies that the map
\[
(\phi_0,\phi_1)\mapsto D^*A\phi(\cdot)\in L^2(\Sigma) 
\]
admits a bounded extension to 
$Y$ for every $T>0$.
In particular, if $(\xi,\eta)\in Y$ then
\begin{equation}
\ZLA{eq:proprieEXRENZ} 
\left\{\begin{array}{l}
  D^*AR_+(t)\xi=\ZT R_+(t)\xi\in L^2(\Sigma)\,,\\
   D^*\ZA R_-(t)\eta=\ZT \ZA^{-1}R_-(t)\eta\in L^2(\Sigma)\,.
   \end{array}\right.
\end{equation}

\subsection{Controllability and Riesz sequences}
Let $ \{\phi_n\} $ be an orthonormal basis of $ L^2(\ZOMq)  $ where $ \phi_n $ is an eigenvector of the operator $A$ with eigenvalue $ \mu_n=\zl_n^2 $. 
We introduce the following sequence of functions $\{\Psi_n\}$:
\begin{equation}
\ZLA{eq:defiPsin}
\left.\begin{array}{ll}
\mbox{\bf case (A)} & \Psi_n=\frac{\ZT\phi_n}{\zl_n^{3/2}}\\[2mm]
\mbox{\bf case (A)}& \Psi_n=\frac{\ZT\phi_n}{\zl_n }\,.
\end{array}\right.
\end{equation}
We note the following inequality which extends the result in~\cite{Tao}:
\begin{Lemma}\ZLA{lemma:almostNORMAL}
 The sequence
$
\{\Psi_n\}$
 is \emph{almost normalized } in $ L^2(\Gamma) $, i.e. there exist $m_0>0$ and $M$ such that
 \[
m_0\leq\|\Psi_n\|_{L^2(\Gamma)}\leq M\,. 
 \] 
\end{Lemma}
The proof is similar to that of~\cite[Lemma~4.4] {PandLIBRO}. We present it  in the case {\bf (A).}
We solve Eq.~(\ref{eq:sistemaelastico}) with zero boundary condition (so the solution is denoted $\phi$), $F=0$ and initial conditions
\[
\phi(0)=\phi_n\,,\quad \phi'(0)=0\,.
\]
The solution is $\phi(x,t)=\left (\cos\zl_n t\right ) \phi_n(x)$ and
\[
\|\phi_n\|_{{\rm dom}\, A^{3/4}}=\zl_n ^{3/2}\,.
\]
Hence, we have
\[
m_0\zl_n^3\leq \left (\intT  \cos^2\zl_n t \ZD t\right )\left (   \int _\Gamma  |\ZT \phi_n(x)|^2 \ZD\Gamma\right )\leq M\zl_n^3
\]
from which the required property follows.
 The proof in the case {\bf (B)} is similar, since $\|\phi_n\|_{{\rm dom }\, A^{1/2}}=\zl_n$.

Now we identify a Riesz sequence which is naturally associated to the controllability of the elastic system~(\ref{eq:sistemaelastico}).

We expand in series of eigenfunctions  the solutions of problem~(\ref{eq:sistemaelastico}), (\ref{eq:elasticoBoundCondi}) (with $ u_0=0 $, $ u_1=0 $) and we find a Fourier representation for the map $ \Lambda_T g $. Let
\[ 
u_n(t)=\int_\ZOMq u(x,t)\phi_n(x)\ZD x\,.
 \]
 Then  (see~\cite{PandLIBRO,PandLame} for similar computations)
\begin{align}
u_n(T)&= -\intT \int_\Gamma g(x,T-s) \frac{1}{\zl_n} \left (\sin \zl_n s\right ) \ZT\phi_n \ZD\Sigma\,,\\
u_n'(T)&= -\intT \int_\Gamma g(x,T-s)\left (\cos\zl_n s\right )\ZT\phi_n\ZD\Sigma\,.
\end{align}

Now we consider case {\bf (A)} so that
   $ \left (u_n(T),u_n'(T)\right )\in \left (\left ({\rm dom}\, A^{1/4}\right )' ,  \left ({\rm dom}\, A^{3/4}\right )' \right )$
and   any $ \xi\in  \left ({\rm dom}\, A^{1/4}\right )'$, $ \eta\in   \left ({\rm dom}\, A^{3/4}\right )' $ has the representation
\[ 
\xi=\ZSUno \sqrt \zl_n \xi_n\phi_n(x)\,,\qquad \eta=\ZSUno \zl_n^{3/2} \eta_n\phi_n(x)\,,\qquad \{\xi_n\}\in l^2\,,\ \{\eta_n\}\in l^2\,.
 \]
So, controllability is equivalent to the surjectivity of the following map $ \Mo $: $ L^2(0,T;L^2(\Gamma))\mapsto l^2\times l^2 $:
 \begin{align*} 
 \Mo g=\left \{\intT \int_\Gamma \left ( \Psi_n  \sin\zl_n s  \right )g(x,T-s)\ZD s\,\ZD\Gamma\,,\right.  \\
\left. \intT \int_\Gamma
 \left (\Psi_n\cos\zl_n s   \right )g(x,T-s)\ZD s\,\ZD\Gamma
 \right\} 
  \end{align*}
  where $\Psi_n$ is given by the first expression in~(\ref{eq:defiPsin}).
  
  If we are in the case~{\bf(B)} then 
  \[
\xi=\ZSUno\xi_n\phi_n(x)\,,\quad \eta=\ZSUno \zl_n\eta_n\phi_n(x)\,,\qquad \{\xi_n\}\in l^2\,,\ \{\eta_n\}\in l^2 
  \]
and  surjectivity of the   operator $\Mo$ (whith $\Psi_n$ given by the second expression in~(\ref{eq:defiPsin}))  is equivalent to controllability also in the case {\bf (B).}

Continuity of the transformation $g\mapsto (u,u')\in C([0,T],X)$ implies continuity of $ \Mo $ so that surjectivity is equivalent to the integration in the definition of $\Mo$ being against a Riesz sequence of functions (see~\cite[Ch.~3]{PandLIBRO}). Hence we have:
  \begin{Theorem}
  The sequence 
  $
    \left \{ \Psi_n e^{i\zl_n t}\right \}
    $
   is a Riesz sequence in $ L^2(0,T;L^2(\Gamma)) $.
  \end{Theorem}

\section{\ZLA{sect:TheSOlutions}The solutions of the viscoelastic system} 
 Let $R(t)$ be the resolvent kernel of $M(t)$, given by  
\begin{equation}\ZLA{eq:resolvM}
R(t)+\intt M(t-s)R(s)\ZD s=M(t)\;.
\end{equation} 
 We ``solve'' the Volterra integral equation~(\ref{eq:sistemsvisco})
 in the ``unknown'' $\Delta^2 w $ (this is called the \emph{MacCamy trick}). 
We get
\[
\Delta^2 w(t)=-w''(t) -\intt R(t-s)  w''(s) \ZD s+F(t)=\intt R(t-s)F(s)\ZD s\,.
\]

 We integrate by parts and we get   
 
\begin{equation}\ZLA{eq:foPoMacCam}
\left\{\begin{array}{l}
\displaystyle    w(t)'' + \Delta^2 w(t)=aw'(t)  +b  w(t)+\intt K(t-s)  w(s)\ZD s + F_1(t)\,,\\
\displaystyle   F_1(t) =  -R(t) w_1-R'(t) w_0+F(t)-\intt R(t-s)F(s)\ZD s \,,\\[2mm]
\displaystyle   a= R(0)\,, \quad 
 b=R'(0)\,,\quad K(t)=R''(t) 
\end{array}\right. 
 \end{equation}
 and the initial and boundary conditions in~(\ref{eq:sistemsvisco}).
 i.e. MacCamy trick removes the differential operator from the memory term.
  
 The term $a w'(t)$ can be removed too from the right hand side if we perform the transformation $ v(t)=e^{-(a/2)t} w(t)$. 
The effect on the initial conditions is   that $ v(0)= v_0= w_0$ while $ v'(0)= v_1= w_1-(a/2) w_0$ and $F_1(t)$, $g(t)$ are replaced by $e^{-(a/2)t}F_1(t)$, $e^{-(a/2)t} g(t)$. Of course this has no influence on controllability and  we assume $a=0$ in~(\ref{eq:foPoMacCam}) from the outset.

We define
\begin{equation}\ZLA{eq:formaL}
L(t)w= \left [
bR_-(t)w+\intt K(t-r)R_-(r)w\ZD r\right ]\,.
\end{equation}
By using~(\ref{eq:soluELASTCO})  we see that $w$ solves also the Volterra integral equation
\begin{equation}
\ZLA{eq:soluzioneviscoelastico}
w(t)=u(t) +\ZA^{-1} \intt L(t-s)w(s)\ZD s
\end{equation}
where
  $u(t)$,  given by~(\ref{eq:soluELASTCO}) with   
  $u_0=w_0$, $u_1=w_1$, $F(t)=F_1(t)$, depends on the boundary control $g$ too.

By definition, \emph{the solutions $w$ of the Volterra integral equation~(\ref{eq:soluzioneviscoelastico}) are the solutions of the Volterra integral equation~(\ref{eq:sistemsvisco}).
}
 
 \begin{Remark}
{\rm
From now on, when convenient for clarity, convolution is also denoted $*$, hence
\[
L*u=\intt L(t-s) u(s)\ZD s\,, \qquad L^{(*k)}*u=L*\left (L^{(*(k-1))}*u\right )\,,\quad k\geq 2\, \zdiaform
\]
}
\end{Remark}
  The operators in~(\ref{eq:soluzioneviscoelastico}) are continuous, and $L(t)$ is strongly continuous, so that we know from the theory of the Volterra integral equations in Hilbert spaces that the regularity properties of $u$ are inherited by $w$:
The regularity properties in Theorem~\ref{teo:proprieREGOela} hold also for $w(t)$.  
  
  Using Picard method, the solution $w(t)$ of~(\ref{eq:soluzioneviscoelastico}) is
  \begin{align}
\nonumber
  w(t)&=u(t)+\ZA^{-1}\intt L(t-s)u(s)\ZD s+A^{-1}\left [\sum _{n=2}^{+\ZIN}\ZA ^{-n+2}\left (L^{(*n)}\right )*u\right ](t)=\\
  \ZLA{eq:SoluVIApicard}  &=u(t)+\ZA^{-1} \intt H(t-s) u(s)\ZD s\,.
  \end{align}

 It is also true that the direct
   inequality~(\ref{eq:direct}) holds for the viscoelastic system:
   \begin{equation}\ZLA{eq:directVISCO}
\int_\Sigma \left |\ZT\psi\right |^2\ZD\Sigma\leq M\left (\left \|(\psi_0,\psi_1)\right \|^2_Y+\|F_1\|^2_{L^1(0,T;L^2(\ZOMq))}
\right ) 
\end{equation}
where $\psi$ solves~(\ref{eq:sistemsvisco}) (with $g=0$).  Of course the constant $M$ is not the same as in~(\ref{eq:direct})). This   easily follows from~(\ref{eq:SoluVIApicard}) that we rewrite  (also for future reference) in terms of $\psi$ and $\phi$:

\begin{align}\nonumber
\psi(t)&=\phi(t)+\ZA^{-1} \intt H(t-s) \phi(s)\ZD s=\\
&=\phi(t)+\ZA^{-1}\intt L(t-s)\phi(s)\ZD s+A^{-1}\left [\sum _{n=2}^{+\ZIN}\ZA ^{-n+2}\left (L^{(*n)}\right )*\phi\right ](t) \,.
\ZLA{eQ:UsOPiCaRd} 
\end{align}

It is sufficient to note that the direct inequality holds for each one of the three addenda. This is clear for the first one (which solves the equation without memory) and for the last 
 term since $\ZT A^{-1}=-D^*$   is bounded, see~(\ref{eq:operatorefrontiera}). 
Now we study
\[
\ZA^{-1}\intt L(t-s)\phi(s)\ZD s\,.
\] 
  Let us note that
\begin{align*}
\phi(t)&=R_+(t)\psi_0+\ZA^{-1}R_-(t)\psi_1-\ZA^{-1}\intt R_-(t-s)F_1(s)\ZD s \\
F_1(s)&= \left [ R(s)\psi_1+R'(s)\psi_0\right ]+ \left [ F-R*F\right ](s) \,.
\end{align*}

The contribution of the last term is
\[
\ZT A^{-1}\intt L(t-\zt)\int_0^\zt R_-(\zt-s) F_1(s)\ZD s\,\ZD\zt
\]
and $\ZT A^{-1}$ is bounded, extended by $-D$. So, we remain with
\begin{align*}
\ZA^{-1} \intt L(t-s)R_+(s)\psi_0(s)\ZD s+A^{-1} \intt L(t-s) R_-(s)\psi_1\ZD s-\\
- A^{-1} \intt L(t-s)\ints R_-(s-\zt)\left [ R(\zt)\psi_1+R'(\zt)\psi_0\right ]\ZD\zt\,\ZD s 
\end{align*}
and in fact only the first addendum: 
\begin{align*}
 \ZA^{-1} \intt L(t-s)R_+(s)\psi_0(s)\ZD s= b\ZA^{-1} \intt R_-(t-s)R_+(s)\psi_0\ZD s+\\
+\ZA^{-1} \intt \int_0^{t-s} K(r)R_-(t-s-r)R_+(s)\psi_0\ZD r\,\ZD s\,.
\end{align*}
We  examine the first integral on the right hand side, which can be computed explicitly:

  \[
 \ZA^{-1} \intt bR_-(t-s) R_+(s)\psi_0\ZD s  =\frac{t}{2} \ZA^{-1}R_-(t)\psi_0+\frac{1}{2}A^{-1}R_-(t)\psi_0\,.
  \] 
 The function $t\mapsto \ZA^{-1} R_-(t)\psi_0$ is the solution $\phi$ of~(\ref{eq:sistemaelastico}) with $g=0$ $F=0$, $\phi_0=0$ and $\phi_1=\psi_0$. Hence (we give the computation in the case {\bf (A).} Case {\bf (B)} is similar)
 \[
 \|\ZT \ZA^{-1}R_-(t)\psi_0\|_{L^2(\Sigma)}\leq M\|\psi_0\|_{{\rm dom} A^{1/4}}< M\|\psi_0\|_{{\rm dom} A^{3/4}} \,.
 \]
 The remaining term is treated analogously.
 
 \emph{This ends the proof of the inequality~(\ref{eq:directVISCO}).}
 \begin{Remark}
 {\rm We shall have the occasion to use~(\ref{eq:foPoMacCam}) with $F_1=0$, which is Eq.~(\ref{eq:sistemsvisco}) with $F$ given by $F-R^*F=R(t)w_1+R'(t) w_0$. Needless to say, the direct inequality holds in this case too.\zdia
 }
 \end{Remark}
\section{\ZLA{sect:ControllabilityVISCO}The proof of controllability}

We assume $F=0$, $u_0=0$, $u_1=0$.

We recall the notation $\Lambda_{T} g=(u_g(T),u_g'(T))$ and we introduce

 \begin{align*}
&  \Lambda_{T}^Vg=(w_g(T),w_g'(T))\,,\\
&R(T)=\left \{ (u_g(T),u_g'(T))\,,\ g\in L^2(0,T;L^2(\Gamma))\right \}\,,\\
& R_V(T)=\left \{ (w_g(T),w_g'(T))\,,\ g\in L^2(0,T;L^2(\Gamma))\right \}\,.
 \end{align*}
The notation is redundant since $R(T)=X$.
We want to prove that, for every $T>0$,
\[
R_V(T)=X\quad {\rm i.e.}\quad 
{\rm im}\, \Lambda_{T}^V=X\,.
\] 
 The proof is in two step. In the first step we prove that $  \Lambda_{T}^V=\Lambda_T+\mathcal{K}_T$ and that $\mathcal{K}_T$ is compact for every $T>0$. We know that $\Lambda_T$ is surjective and so the image of $ \Lambda_{T}^V$ is closed with finite codimension. Then we prove surjectivity.

\paragraph{Step~1: the codimension of $R_V(T)$ is finite}   
We fix $g\in L^2(\Sigma)$ and we consider 
 the first component $w_g(T)$ of $ \Lambda_{T}^V$.  
We use~(\ref{eq:soluzioneviscoelastico}) 
 
\begin{align*}\nonumber
 w_g(T)&=u_g(t)+\ZA^{-1}\intT H(T-s) u_g(s)\ZD s=\\
 &
 =\Lambda_T g+
 \ZA^{-1}\intT H(T-s) u_g(s)\ZD s\,.
\end{align*}
Now we give the details in the case {\bf (A)} (easily adapted to the case {\bf(B)}).

 Theorem~\ref{teo:proprieREGOela} shows that $u_g\in C\left ([0,T],\left ({\rm dom}A^{1/4}\right) '\right )$ and so
    $g\mapsto (H*u_g)(T)$ is continuous from $L^2(\Sigma)$ to $ \left ({\rm dom}A^{1/4}\right) '$. The operator $\ZA^{-1}$ (extended to $ \left ({\rm dom}A^{1/4}\right) '$ by duality) is compact and so $g\mapsto w_g(T)$ is the sum of a surjective and a compact operator.

  The velocity component $w_g(T)$ is treated analogously by using
  \begin{equation}
  \ZLA{eq:dellaVeloPerRT}
  w'_g(T)=u'_g(T)+\ZA^{-1}\intT H(s)u'_g(T-s)\ZD s 
  \end{equation}
  (here we used $H(0)=0$, as seen from~(\ref{eq:formaL}) and~(\ref{eq:SoluVIApicard})).
   We get that $\Lambda^V_T$ is a compact perturbation of the surjective operator $\Lambda_T$. Hence, the codimension of $R_V(T)$ is finite.
  
  Now we prove $\left  [R_V(T)\right ]^\perp=0$. This part of the proof requires several steps.

 \paragraph{Step 2: characterization of the elements of $\left[R_V(T)\right ]^\perp$}
 We recall $X=\tilde Y'$.
  We characterize the annihilators  in $\tilde Y$   of $R_V(T)$ respect to the duality pairing. The set of the annihilators is denoted   $\left [R_V(T)\right ]^\perp$.

 The pairing $\ZLL\cdot,\cdot\ZRR$ of $\tilde Y$ and its dual $X$  is the inner product in $L^2$ if it happens that $(w(T),w'(T))$ belongs to $L^2(\ZOMq)\times L^2(\ZOMq)$ and this is the case if $g\in \mathcal{D}(\Sigma)$.
  The reachable set with   $g\in \mathcal{D}(\Sigma)$ is dense in $R_V(T)$ and so we confine ourselves to compute with  smooth $g$. In this case,
  
\begin{equation}\ZLA{eq:duality}
\ZLL (\eta,\xi),(w(T),w'(T))\ZRR=\int_\ZOMq \eta w(T)\ZD x+\int_\ZOMq \xi w(T)\ZD x\,.
\end{equation}
We examine the first integral, using~(\ref{eq:SoluVIApicard})
\begin{align*}
&-\int_\ZOMq\eta w(T)\ZD x=\int_\ZOMq \eta(x)\ZA \intT R_-(T-s)Dg(s)\ZD s\,\ZD x+\\
&+ \int_\ZOMq \eta \intT \int_0^{T-s}  R_-(T-s-r)H(r) Dg(s)\ZD r\,\ZD s\,\ZD x=\\
&= \int_\Sigma g(s) D^*\ZA R_-(T-s)\eta\ZD\Sigma+\int _\Sigma g(s)D^* \int_0^{T-s}R_-(T-s-r) H(r)\eta \ZD r\,\ZD\Sigma=\\
&=\int_\Sigma g(T-s)D^*A\left [
\ZA^{-1} R_-(s)\eta+A^{-1} \ints R_-(s-r)H(r)\eta\ZD r
\right ]\ZD\Sigma\,.
\end{align*}
Note that the previous computations are justified by~(\ref{eq:proprieEXRENZ}).

The function
\begin{equation}
\ZLA{eq:primaDIpsi}
\psi(s)=\ZA^{-1 } R_-(s)\eta+\ZA^{-1} \ints H(s-r)\ZA^{-1} R_-(r)\eta \ZD r
\end{equation}
is~(\ref{eQ:UsOPiCaRd}) with $\phi(s)=\ZA^{-1}R_-(s)\eta$. Hence $\psi(s)$ it is the solution   of~(\ref{eq:foPoMacCam}) (where $a=0$) with $F_1=0$, $g=0$, $w(0)=\psi(0)=0$, $w'(0)=\psi'(0)=\eta$ 
and we have
 \begin{equation}\ZLA{eq:OrtoPriCONDI}
-\int_\ZOMq \eta w(T)\ZD x=\int_\Sigma g\ZT \psi\ZD\Sigma \,.
 \end{equation}

 The integral 
 $\int_\ZOMq \xi w'(T)\ZD x$ is treated analogously. We use~(\ref{eq:dellaVeloPerRT}). We see that
\[
\ZLA{expreDERIVAprima}
w'(t)=-A\intt R_+(t-s)Dg(s)\ZD s-\ZA\intt R_+(t-s)\ints H(s-r) Dg(r)\ZD r\,\ZD s 
\]
and we get
 \begin{align*}
 &-\int_\ZOMq \xi w'(T)\ZD x=\int_\ZOMq \xi A\intT R_+(T-s)Dg(s)\ZD s\,\ZD x+\\
 &+\int_\ZOMq \xi\ZA\intT R_+(T-s)\ints H(s-r)Dg(r)\ZD r\ZD s\,\ZD x\,.
 \end{align*}
 As above, we use~(\ref{eq:proprieEXRENZ}) in order to justify the following exchange of integration:
 \begin{align}
 \nonumber  &-\int_\ZOMq \xi w'(T)\ZD x=\\
 \nonumber  &= \int_\Sigma g(T-s)D^*A\left [ R_+(s)\xi+\ZA^{-1} \ints H(s-r) R_+(r)\xi\ZD r\right ]\ZD\Sigma=\\
 \ZLA{eq:OrtoSecoCONDI}  &= \int_\Sigma g\ZT \psi\ZD\Sigma
 \end{align}
 where now
 \[
\psi(s)=  R_+(s)\xi+\ZA^{-1} \ints H(s-r) R_+(r)\xi\ZD r
 \]
is~(\ref{eQ:UsOPiCaRd}) with $\psi(s)= R_+(s)\xi$. Hence $\psi$ solves~(\ref{eq:foPoMacCam}) (where $a=0$) with $F_1=0$, $g=0$, $w(0)=\psi(0)=\xi$, $w'(0)=\psi'(0)=0$.
 
If $(\eta,\xi)\in \left [R_V(T)\right ]^\perp$ then the sum of (\ref{eq:OrtoPriCONDI}) and (\ref{eq:OrtoSecoCONDI}) is zero for every $g$. So, we have
 \begin{Theorem}
 The pair $(\xi,\eta)\in \tilde Y$ annihilates $R_V(T)$ if and only if the solution $\psi$ of
 \begin{equation}\ZLA{eq:aggiunta}
\left \{
\begin{array}
{l}
\psi''+\Delta^2\psi=b\psi+\intt K(t-s)\psi(s)\ZD s\\
\psi(0)=\xi\,,\quad \psi'(0)=\eta\\
\zg_0\psi(t)=0\,,\qquad  \zg_1\psi=0
\end{array}
\right .
\end{equation}
 satisfy the additional condition
\begin{equation}
\ZLA{eq:CondiORTOG}
 \ZT\psi =0\quad {\rm in}\ L^2(\Sigma))\,.
\end{equation}
 
 \end{Theorem}

The goal now is the proof that any solution of~(\ref{eq:aggiunta}) which satisfy~(\ref{eq:CondiORTOG}) is identically zero.

\paragraph{Step 3: regularity of the elements of $\left [R_V(T)\right ]^\perp$}

Here the special definition of the space  $X$ in the cases {\bf (A)} and~{\bf (B)} has a role. We consider the case  {\bf (A)}  first. 

We expand $\xi$ and $\eta $ in series of the eigenfunctions of $A$. The conditions $\xi\in {\rm dom}\, A^{3/4}$, $\eta\in {\rm dom}\, A^{1/4}$ give
 
 \begin{equation}\ZLA{eq:FourieEXPA}
\xi=\ZSUno  \xi_n\phi_n=\ZSUno \frac{\tilde\xi_n}{\zl_n^{3/2}}\phi_n\,,\qquad \eta=\ZSUno \eta_n\phi_n=\ZSUno \frac{\tilde\eta_n}{\zl_n^{1/2}}
 \end{equation}
 and 
 \begin{equation}
 \ZLA{eq:gliXIedETAnelCasoA}
\left \{\tilde \xi_n \right \}=\left \{\zl_n^{3/2}\xi_n\right \} \in l^2\,,\  \left \{\tilde \eta_n\right \}=\left \{\zl_n^{1/2}\eta_n\right \}\in l^2 \,.
 \end{equation}

We expand the solution $\psi$ of~(\ref{eq:aggiunta}):
\[
\psi(t)=\ZSUno \psi_n(t)\phi_n(x)\ \ {\rm where}\  \psi_n''=-\zl_n^2\psi_n+b\psi_n+\intt K(t-s)\psi_n(s)\ZD s 
\]
so that
\[
\psi_n'(t)=-\zl_n^2\intt\psi_n(s)\ZD s+\intt \mathcal{H}(t-r)\psi_n(r)\ZD r+\eta_n\,, \qquad \psi_n(0)=\xi_n 
\]
where 
\[
\mathcal{H}(t)=b+\intt K(s)\ZD s
\]
(note that $\mathcal{H}$, not to be confused with $H$ in~(\ref{eQ:UsOPiCaRd}), is a real valued function).
We introduce the solution $z_n(t)$ of
 
i.e.
\[
z_n'(t)=-\zl_n^2\intt z_n(s)\ZD s+\intt \mathcal{H}(t-t) z_n(r)\ZD r\,,\qquad z_n(0)=1 
\]
i.e.
\[
z_n''=-\zl_n^2 z_n +b z_n+\intt K(t-s) z_n(s)\ZD s\,,\qquad z_n(0)=1\,,\ z_n'(0)=0\,.
\]
Then we have
 
\[
 \psi_n(t)=\xi_n z_n(t)+\left [\intt z_n(s)\ZD s\right ]\eta_n 
\]
and
\[
\psi(t)=\ZSUno\phi_n\left \{ \xi_n z_n(t)+\eta_n\left [ \intt z_n(s)\ZD s\right ]\right \}\,.
\]
The direct inequality justifies the exchange of $\ZT$ and the series, so that the condition of orthogonality is
\begin{align}\nonumber
\ZT\psi&=\ZSUno \left (\ZT\phi_n\right )\left [ \xi_n z_n(t)+\eta_n\left [ \intt z_n(s)\ZD s\right ]\right ]=\\
\ZLA{eq:condiORTOGserie}
&=
\ZSUno
\Psi_n
\left [ \tilde \xi_n z_n(t)+\tilde \eta_n\left (\zl_n\intt z_n(s)\ZD s\right )\right ]=0
\,.
\end{align}
With the proper definition  of $\Psi_n$ (in~(\ref{eq:defiPsin})) this is also the orthogonality condition in the case~{\bf (B)} since in this case
\begin{align}
\nonumber&\xi=\ZSUno\xi_n\phi_n(x)=\ZSUno\frac{\tilde\xi_n}{\zl_n} \phi_n(x)\,,\  \eta=\ZSUno \eta_n\phi_n(x)\,,\\
\ZLA{eq:gliXIedETAnelCasoB}&\{\tilde\xi_n\}=\{\zl_n\xi_n\}\in l^2\,,\ \{\tilde\eta_n\}=\{\eta_n\}\in l^2\,.
\end{align}

The sequences  $\{\tilde \xi_n\}$ and $\{\tilde\eta_n\}$ in~(\ref{eq:condiORTOGserie}) belong to $l^2$ but in fact they are more ``regular''.   
To see this, we use the following representation of $z_n(t)$:
\begin{align}
\nonumber \displaystyle &z_n(t) = \cos\zl_n t+\frac{b}{\zl_n} \intt \sin\zl_n(t-s) z_n(s)\ZD s+\\ 
\nonumber\displaystyle &+\frac{1}{\zl_n}\intt \left [\int_0^{t-s} \sin\zl_n(t-s-r) K(r)\ZD r\right ]z_n(s)\ZD s\,,\\ 
\nonumber\displaystyle &\intt z_n(\zt)\ZD \zt =  -\frac{1}{\zl_n} \sin\zl_n t-\frac{b}{\zl^2_n} \intt \left [1-\cos\zl_n(t-s)\right ] z_n(s)\ZD s-\\ 
\ZLA{eq:RAppreAsintoZ} \displaystyle &- \frac{1}{\zl_n^2}\intt \left [1-\cos\zl_n(t-s)\right ] \ints K(s-r)z_n(r)\ZD r\,\ZD s\,.
\end{align}
So, the condition of orthogonality is
 
\begin{align}
\nonumber  &\ZSUno \Psi_n
 \left [\tilde\xi_n\cos\zl_n t-\tilde \eta_n\sin\zl_n t\right ]=\\
\nonumber   & =-\ZSUno  \Psi_n\left \{
  \tilde\xi_n \frac{b}{\zl_n}\intt 
  \sin\zl_n(t-s) z_n(s)\ZD s \right.+\\
\nonumber  &\left.+\frac{\tilde\xi_n}{\zl_n}\intt \left [\int_0^{t-s} \sin\zl_n(t-s-r)K(r)\ZD r\right    ] z_n(s)\ZD s\right\}+\\
\nonumber   & +\ZSUno 	 \Psi_n
 \left \{\tilde\eta_n \frac{b}{\zl_n}\intt    \left [1-\cos\zl_n(t-s)\right ]z_n(s)\ZD s+\right. \\
\ZLA{eq:formaPERc1} &\left.+\tilde \eta_n\frac{1}{\zl_n}\intt  \left  (1-\cos\zl_n(t-s)\right ) \ints   K(s-r)z_n(r)\ZD r\right \}   
\,.
\end{align}

Now we prove:
\begin{Lemma}
 The right hand side of~(\ref{eq:formaPERc1})  is of class $H^1$ and so there exist $l^2$ sequences
 $\left \{\widetilde{\tilde \xi }_n\right \}$,  $\left \{\widetilde{\tilde \eta }_n\right \}$  such that 
 \begin{equation}\ZLA{eq:FormaDIxiedETA}
\tilde \xi_n= \frac{1}{\zl_n}\widetilde{\tilde \xi }_n\,,\qquad 
\tilde \eta_n= \frac{1}{\zl_n}\widetilde{\tilde \eta }_n\,,
 \end{equation}
\end{Lemma}  
\zProof
In this proof we use explicitly $\{1/\zl_n\}\in l^2$ so that $\{\tilde\xi_n/\zl_n\}$ and $\{\tilde\eta_n/\zl_n\}$ belong to $l^1$. This condition holds since ${\rm dim}\,\ZOMq=2$, see~(\ref{eq:stimaautov}) but, as we noted, it can be removed.

We distribute on the sum  the series  on the right hand side. This is legitimate since Lemma~\ref{lemma:almostNORMAL} and~(\ref{eq:stimaautov}) show convergence of the resulting series (we use also that $\{z_n(t)\}$ is bounded on bounded intervals).

We compute termwise the derivatives of the individual series  and we prove convergence of the series of the derivatives.
We show the method on the first series (the other ones can be treated similarly).

The derivative   is
\begin{align*}
&b\ZSUno \Psi_n\tilde \xi_n\intt\cos\zl_n(t-s) z_n(s)\ZD s=\\
&b\ZSUno  \Psi_n\left \{ \tilde\xi_n\intt \cos\zl_n(t-s)\cos\zl_n s\ZD s+  \frac{M_n(t)}{\zl_n}\tilde\xi_n\right \}
\end{align*}
and $\{M_n(t)\}$ is bounded. So,    $\ZSUno \Psi_n \frac{M_n(t)}{\zl_n}\tilde\xi_n$ is uniformly convergent.

We note that
\begin{align*}
\ZSUno  \Psi_n \tilde\xi_n\intt \cos\zl_n(t-s)\cos\zl_n s\ZD s=\\
=\left ((1/2)t\right )\ZSUno \Psi_n  \tilde\xi_n  \cos\zl_n t+\ZSUno  \Psi_n \frac{1}{\zl_n} \tilde\xi_n\sin\zl_n t\,.
\end{align*}
The last series is uniformly convergent while the first series on the right hand side converges in $L^2(0,T)$ since
 $\left  \{ \Psi_n e^{i\zl_n t }\right \}$ is a Riesz sequence in $L^2(0,T )$.

The remaining terms  are treated analogously.

In conclusion the right, hence also the left hand side of~(\ref{eq:formaPERc1}) is of class $H^1(0,T)$. We recall  that $\left \{ \Psi_ne^{-\zl_n t}\right \}$ is a Riesz sequence in $L^2(0,T-\ZEP)$ for every $T>0$ and   $\ZEP\in (0,T)$ (since the elastic system is controllable for every $T>0$) and the conclusion follows from~\cite[Lemma~3.4]{PandLIBRO}.\zdia

Note that in particular this result shows uniform convergence of the series~(\ref{eq:formaPERc1}) and so, computing with $t=0$ we get 
\[
\ZSUno \xi_n\ZT\phi_n    =0\,.
\]
Now we replace in the derivative of~(\ref{eq:formaPERc1}) the expressions~(\ref{eq:FormaDIxiedETA}) of $\tilde\xi_n$ and $\tilde\eta_n$  and we see that we can do a second derivative. In conclusion we get:
\begin{Theorem}\ZLA{Teo:regolarDELpasso3}
There exist $l^2$ sequences
 $\left \{  \hat \xi  _n\right \}$,  $\left \{  \hat \eta  _n\right \}$  such that 
 \[
\tilde\xi_n= \frac{1}{\zl_n^2} {\hat \xi }_n\,,\qquad 
\tilde \eta_n= \frac{1}{\zl_n^2} {\hat \eta }_n\,.
 \]
 Furthermore we have
 \begin{equation}\ZLA{eq:UguaAzero}
\ZSUno \ZT  \phi_n \xi_n  =0\,,\qquad \ZSUno  \ZT\phi_n     \eta_n=0\,.
 \end{equation}
 \end{Theorem}

\paragraph{Step 4: More elements in $\left [R_V(T)\right ]^\perp$}
 
The regularity results we found in the previous step allows termwise computation of the derivative of
 $\ZT\psi(t)$ which of course   is zero. We do the computation using formula~(\ref{eq:condiORTOGserie}) and we use the notation $\Psi_n=(1/\zl_n^{3/2})\mathcal{T}\phi_n$. We get:
 \[
\ZSUno\Psi_n \left [\tilde\xi_n z_n'(t)+ \tilde\eta_n\zl_n+\tilde\eta_n\zl_n\intt z_n'(t-s)\ZD s\right ]=0\,.
 \]
 The last equality in~(\ref{eq:UguaAzero}) shows that $\ZSUno \Psi_n \tilde \eta_n\zl_n=0$. So, from the definition of $z_n(t)$ we see that:
 
\begin{align}
\nonumber  &   \ZSUno \Psi_n \left [\tilde\xi_n
\left (
-\zl_n^2\intt z_n(s)\ZD s+\intt \mathcal{H}(t-s)z_n(s)\ZD s
\right )+\right.\\
\nonumber   & +\left.
 \tilde\eta_n\zl_n\left (
\intt \left (
-\zl_n^2\int_0^{t-s} z_n(r)\ZD r+\int_0^{t-s}\mathcal{H}(t-s)z_n(r)\ZD r
\right )\ZD s
\right )
\right ]=\\
\nonumber  &  =\intt \left [  \ZSUno \Psi_n  \left (   -\zl_n^2\tilde\xi_nz_n(s)+\left (-\zl_n^2\tilde \eta_n  \right )                 \zl_n \ints z_n(r)\ZD r \right )               \right ]\ZD s+\\
\nonumber  & +\intt \mathcal{H}(t-s)\left [ 
\ZSUno\Psi_n\left (
\tilde \xi_n z_n(s)+\tilde\eta_n\zl_n\ints z_n(r)\ZD r
\right )
\right ]\ZD s=\\
\ZLA{eq:priMAparteNUOveELEMENTO}  &= -\intt \left [  \ZSUno \Psi_n  \left (   \hat\xi_nz_n(s)+   \hat\eta_n                   \zl_n \ints z_n(r)\ZD r \right )               \right ]\ZD s=0
\end{align}
The series can be distributed thanks to~(\ref{eq:RAppreAsintoZ})  since   $\left \{\Psi_n e^{i\zl_n t}\right \} $ is a Riesz sequence and $\{\zl_n^2\tilde\xi_n\}=\{\hat\xi_n\}\in l^2$, $\{\zl_n^2\tilde\eta_n\}=\{\hat\eta_n\}\in l^2$.

The equality at the line~(\ref{eq:priMAparteNUOveELEMENTO}) is the same as~(\ref{eq:condiORTOGserie}) but for the sequences $\{\hat\xi_n\}$, $\{\hat\eta_n\}$. So, we get a second element $(\xi_1,\eta_1)\in \left [R_V(T)\right ]^\perp$, given by (compare~(\ref{eq:FourieEXPA}) and~(\ref{eq:gliXIedETAnelCasoB}))
\[
\begin{array}{lll}
\mbox{\bf case (A)}&\xi_1=\ZSUno \frac{1}{\zl_n^{3/2}}\hat\xi_n\phi_n(x)\,,& \eta_1=\ZSUno \frac{1}{\sqrt{\zl_n}}\hat\eta_n\phi_n(x)\,,\\[2mm]
\mbox{\bf case (B)}&\xi_1=\ZSUno \frac{1}{\zl_n }\hat\xi_n\phi_n(x)\,,&  \eta_1=\ZSUno\hat\eta_n\phi_n(x)
\end{array}
\]
 We have two cases: either the nonzero elements in both the series correspond to the same eigenvalue, or $(\xi_1,\eta_1)$ (if nonzero) is not colinear with $(\xi,\eta)$. We shall see that the first case cannot happens and so we have constructed \emph{a second element of $\left [R_V(T)\right ]^\perp$ such that  $(\xi,\eta)$, $(\xi_1,\eta_1)$ are linearly independent.}
 
 We repeat  the same computations with $(\xi,\eta)$ in~(\ref{eq:condiORTOGserie}) replaced by $(\xi_1,\eta_1)$.  If the nonzero elements of the series correspond to \emph{at least two different eigenvalues,} we get a third element $(\xi_2,\eta_2)\in\left [R_V(T)\right ]^\perp$  which is linearly independent of the previous ones.
 
 The procedure can be iterated, but only finitely many times, since we already proved that 
 $\left [R_V(T)\right ]^\perp$ is finitely dimensional.
 
  \emph{This proves   that both  the series~(\ref{eq:FourieEXPA}) are finite sums:}
  if $0\neq (\xi,\eta)\in \left [R_V(T)\right ]^\perp $  then (recall $\zl_n^2\leq \zl_{n+1}^2$)
 \begin{equation}\ZLA{eq:FourieEXPAfinite}
\xi=\ZSUno  \xi_n\phi_n=\sum_{\zl_n\leq\zl_N} \xi_n\phi_n\,,\qquad \eta=\ZSUno  \eta_n \phi_n=\sum_{\zl_n\leq\zl_N}  \eta_n\phi_n  
 \end{equation}
 and the orthogonality condition~(\ref{eq:condiORTOGserie}) is
\begin{equation}\ZLA{eq:LaCondortoFIniTA}
 \ZT\psi =\sum_{\zl_n\leq\zl_N} \left (\ZT\phi_n\right )\left [ \xi_n z_n(t)+\eta_n\left ( \intt z_n(s)\ZD s\right )\right ]=0\,.
 \end{equation}

 In the last step we are going to see that the series~(\ref{eq:FourieEXPAfinite}) have to be equal zero.

\paragraph{Step 5: end of the proof}
We recall our goal: if $(\xi,\eta)\in \left [R_V(T)\right ]^\perp$ then we must prove $\xi=0$, $\eta=0$ and we know that   $\xi$ and $\eta$ are given by~(\ref{eq:FourieEXPAfinite}) while the orthogonality condition is~(\ref{eq:LaCondortoFIniTA}). Now the proofs goes as in~\cite{PandParma,PandLame} and is reported for completeness.

We prove that if~(\ref{eq:FourieEXPAfinite})-(\ref{eq:LaCondortoFIniTA}) hold  (with $\|\xi\|^2+\|\eta\|^2\neq 0$) then it is possible to find a new nonzero element of $\left  [R_V(T)\right ]^\perp $ whose nonzero   coefficients  $\xi_n$ and $\eta_n$ all correspond to the one and the same eigenvalue. Let us assume that in the sums~(\ref{eq:FourieEXPAfinite}) we can find nonzero coefficients which correspond to two different eigenvalues $\zl_N^2$ and $\zl_{k_0}^2$.
We compute the derivative of $ \ZT\psi(t)$ (which is zero). The same computation as in the step~4 shows that
\begin{align} 
 \nonumber \left ( \ZT\psi\right )' =\sum_{\zl_n\leq\zl_N} \left (\ZT\phi_n\right )\left [ \zl_n^2\xi_n z_n(t)+\zl_n^2\eta_n\left ( \intt z_n(s)\ZD s\right )\right ]=\\
\ZLA{eq:LaCondortoFIniTAMAlambdaQUADRO}=\sum_{\zl_n=\zl_N} +\sum_{\zl_n=\zl_{k_0}} +\sum _{\zl_n<\zl_N\,, \zl_n\neq\zl_{k_0}}=0
 \end{align} 
Note that ($m_N$ is the multiplicity of $\zl_N$)
\[
\sum_{\zl_n=\zl_N}=\zl_N^2z_N(t)\left (\sum _{k=1}^{m_N} \left (\ZT\phi_k\right ) \xi_k\right )+ \zl_N^2\left (\intt z_N(s)\ZD s\right )\sum _{k=1}^{m_N} \left (\ZT\phi_k\right ) \eta_k 
\]
 and that $\ZT\psi$ has similar expressions, without the coefficients $\zl_n^2$.
 
 A linear combination of~(\ref{eq:LaCondortoFIniTA}) and~(\ref{eq:LaCondortoFIniTAMAlambdaQUADRO}) with coefficients respectively $\zl_N^2$ and $-1$ cancel the first sum in the right hand side of~(\ref{eq:LaCondortoFIniTAMAlambdaQUADRO}) and we get
 \begin{align*}
 \sum_{\zl_n=\zl_{k_0}} \left (\ZT\phi_n\right )\left [ (\zl_N^2-\zl_{k_0}^2)\xi_n z_n(t)+(\zl_N^2-\zl_n^2)\eta_n\left ( \intt z_n(s)\ZD s\right )\right ]+
\\
 +\sum_{\stackrel{\zl_n\leq\zl_N-1}{ \zl_n\neq \zl_{k_0}}} \left (\ZT\phi_n\right )\left [   (\zl_N^2-\zl_{n}^2)  \xi_n  z_n(t)+ (\zl_N^2-\zl_{n}^2) \eta _n\left ( \intt z_n(s)\ZD s\right )\right ] 
=\\
\sum_{\zl_n<\zl_N} \left (\ZT\phi_n\right )\left [    \xi_n^{(1)} z_n(t)+   \eta^{(1)}_n\left ( \intt z_n(s)\ZD s\right )\right ] 
  =0\,.
 \end{align*}
 It follows that the element whose expansion in egenfunctions has the coefficients $ \xi_n^{(1)}$ and $ \eta_n^{(1)}$ (zero when $\zl_n\geq \zl_N$) belongs to $\left  [R_V(T)\right ]^{\perp} $ and it is nonzero because the coefficients which correspond to the eigenvalue $\zl_{k_0} $ are not all equal to zero.

This procedure can be repeated till we get a sum with terms which correspond to only one eigenvalue \emph{and which has  nonzero coefficients.} 
 We prove that this cannot be, and so our original $(\xi,\eta)\in\left  [R_V(T)\right ]^{\perp}  $
is in fact $\xi=0$, $\eta=0$, as we wanted to prove.

Assume that the sums~(\ref{eq:FourieEXPAfinite}), (\ref{eq:LaCondortoFIniTA}) are for $\zl_n=\zl_N$. In this case $\xi$ and $\eta$ are eigenvectors of the operator $A$ which correspond to the same eigenvalue $\zl_N$ and computing~(\ref{eq:LaCondortoFIniTA}) and its derivative with $t=0$ we get
\[
\ZT\xi=\sum _{\zl_n=\zl_N} \xi_n\mathcal{T}\phi_n=0\,,\qquad \ZT\eta= \sum _{\zl_n=\zl_N} \eta_n\mathcal{T}\phi_n=0\,.
\]

Now we solve~(\ref{eq:sistemaelastico}) with $F=0$, $g=0$ (the solution is denoted $\phi$ because $g=0$) and initial conditions 
\[
 \phi(0)=\xi\,,\qquad  \phi'(0)=\eta\,.
\]
The solution is
\[
\phi(t)=\xi\cos\zl_N t+\eta\frac{1}{\zl_N}\sin\zl_N t\quad\mbox{so that}\quad \mathcal{T}\phi(t)=0\,.
\]
The inverse inequality~(\ref{eq:inverse}) of the \emph{elastic system} implies $\xi=0$ and $\eta=0$, as we wanted to achieve.

\emph{The proof is now finished.}

\end{document}